\documentclass[11pt,a4paper]{amsart}
\usepackage{amsmath,amssymb, amsbsy}
\usepackage{subfigure}
\usepackage{graphpap,latexsym,epsf}
\usepackage{color,psfrag}
\usepackage[dvips]{graphicx}
\usepackage{bm}

\usepackage[english]{babel}
\usepackage{times}
\usepackage{graphicx}
\usepackage{amscd}
\usepackage{amsmath}
\usepackage{amsfonts}
\usepackage{amssymb}
\usepackage{amsthm}
\usepackage{latexsym}

\numberwithin{equation}{section}

\newtheorem{theorem}{Theorem}[section]
\newtheorem{proposition}[theorem]{Proposition}

\newtheorem{remark}{Remark}[section]

\def\neweq#1{\begin{equation}\label{#1}}
\def\endeq{\end{equation}}

\newcommand{\N}{\mathbb{N}}

\newcommand{\SI}{{\bm \sigma}}
\newcommand{\e}{{\bf e}}

\begin{document}

\title{A note on an orthotropic plate model describing the deck of a bridge}

\author{Alberto Ferrero}

\address{\hbox{\parbox{5.7in}{\medskip\noindent{Alberto Ferrero, \\
Universit\`a del Piemonte Orientale, \\
        Dipartimento di Scienze e Innovazione Tecnologica, \\
        Viale Teresa Michel 11, 15121 Alessandria, Italy. \\[5pt]
        \em{E-mail address: }{\tt alberto.ferrero@uniupo.it}}}}}

\date{\today}

\maketitle

\begin{abstract} The purpose of this work is to develop a model for a rectangular plate made of an orthotropic material.
If compared with the classical model of the isotropic plate, the relaxed condition of orthotropy increases the degrees of freedom as a consequence of the larger number of elastic parameters, thus allowing to better describe rectangular plates having different behaviors in the two directions parallel to the edges of the rectangle.

We have in mind structures like decks of bridges where the rigidity in the direction of their length does not necessarily coincide with the one in the direction of its width.

We introduce some basic notions from the theory of linear elasticity, having a special attention for the theory of orthotropic materials. In particular we recall the Hooke's law in its general setting and we explain how it can be simplified under the orthotropy assumption.

Following the approach of the Kirchhoff-Love model, we obtain the bending energy of an orthotropic plate and from it the corresponding equilibrium equation when the plate is subject to the action of a vertical load.
Accordingly, we write the kinetic energy of the plate which combined with the bending energy gives the complete Lagrangian; classical variational methods then produces the equation of motion.

\end{abstract}

\vspace{11pt}


\vspace{6pt}

\section{Introduction} \label{s:introduction}

Let us consider a rectangular plate of length $L$, width $2\ell$ and thickness $d$. Denoting by $E$ and $\nu$ the Young modulus and Poisson ratio respectively, and by $u=u(x,y)$ the vertical displacement of the middle surface of the plate, the resulting elastic bending energy reads
\begin{equation}\label{eq:plate-isotropic}
  \mathbb E_B(u)=\frac{Ed^3}{24(1-\nu^2)} \int_{(0,L)\times (-\ell,\ell)}\left(\nu |\Delta u|^2+(1-\nu)|D^2 u|^2\right)dxdy
\end{equation}
where $D^2 u$ denotes the Hessian matrix of $u$ and $|D^2 u|=\left(u_{xx}^2+2u_{xy}^2+u_{yy}^2\right)^{1/2}$ is its Euclidean norm, obtained by interpreting the matrix as a vector of four components. For more details on the plate model see \cite{BeBuoGa, BeFeGa, FeGa, LaLi}.
This is known as the Kirchhoff-Love model for a plate of an isotropic material, see \cite{Kirchhoff, Love}.

Denoting by $f$ a vertical load per unit of surface, the corresponding Euler-Lagrange equation becomes
\begin{equation} \label{eq:plate-iso-EuLa}
\frac{Ed^3}{12(1-\nu^2)} \, \Delta^2 u=f \qquad \text{in } \Omega=(0,L)\times (-\ell,\ell) \, .
\end{equation}
The constant $\frac{Ed^3}{12(1-\nu^2)}$ in front of the biharmonic operator $\Delta^2$ represents the rigidity of the plate.

For more details on recent literature about rectangular plates and applications in models for decks of bridges, we quote \cite{AnGa, BeBuGaZu, BeFaFeGa, BoGaMo, ChGaGa, Suspension, Gazzola-Libro} and the references therein.

As explained in more details in \cite{Suspension}, an isotropic plate could be not optimal for describing the deck of a bridge since the rigidity constant $\frac{Ed^3}{12(1-\nu^2)}$ determines simultaneously the rigidity in both length and width directions.

For this reason, it appears more reasonable to describe the deck of a bridge combining the equation of beams for vertical displacements and the equation of rods for torsion, see \cite{ArLaVaMa, AG15, AG17}.

However, a possible alternative approach could be realized interpreting the deck as an orthotropic plate. This is the approach we want to follow in the present article and that was considered in \cite{Suspension}.

We point out that the present paper can be seen as a sort of auxiliary article to \cite{Suspension} which has the purpose of clarifying in more details some aspects that were treated there and to collect some basic notions well known in the theory elasticity.

As anticipated at the beginning of this introduction, a thin plate can be interpreted as a three-dimensional solid body which in a suitable coordinate system can be represented by the open set $(0,L)\times (-\ell,\ell)\times \left(-\frac d2, \frac d2\right)$.

Denoting by $E_1>E_2=E_3$ the \textit{Young moduli} in the $x$, $y$, $z$ directions respectively, and $\nu_{12}$ and $\nu_{21}$ the \textit{Poisson ratios} relative to the $x$ and $y$ directions, see Section \ref{ss:Ortho-Mat} for the correct definitions, we will show that the equilibrium equation for the orthotropic plate subject to a vertical load $f$ is given by
\begin{equation} \label{eq:familiar}
 \frac{E_1 d^3}{12(1-\nu_{12} \nu_{21})} \, \frac{\partial^4 u}{\partial x^4}
 +\frac{E_2 d^3}{6(1-\nu_{12} \nu_{21})} \, \frac{\partial^4 u}{\partial x^2\partial y^2}
 +\frac{E_2 d^3}{12(1-\nu_{12} \nu_{21})} \, \frac{\partial^4 u}{\partial y^4}=f \, ,
\end{equation}
see also Remark \ref{r:1}.

Let us denote by $\mathcal A$ the linear fourth order operator appearing in the left hand side of \eqref{eq:familiar}.

Section \ref{s:orthotropic} is devoted to a survey of some well known notions about elastic anisotropic materials including the anisotropic Hooke's law and the related {\it stiffness matrix}. We explain how the stiffness matrix looks like in the case of an orthotropic material, we show how the elastic energy per unit of volume can be expressed in terms of the components of the strain tensor and we describe the meaning of the elastic coefficients involved in the orthotropic Hooke's law, still named Young moduli (one for each coordinate axis) and Poisson ratios (one for each combination of two coordinate axes); other coefficients completely independent of the previous ones are the so called {\it shear moduli}.
We complete that section with a description of an orthotropic material with a two-dimensional symmetry (with respect to planes parallel to the $yz$ plane) and a one-dimensional reinforcement in the orthogonal direction (the $x$ axis).

Having in mind the explanations given in Section \ref{s:orthotropic} about orthotropic materials, in Section \ref{s:orthotropic-plate} we construct a model of orthotropic plate with a one-dimensional reinforcement giving rise to the equation of the orthotropic plate. Since our purpose is to describe the deck of a bridge, the rectangular plate is supposed to be hinged at the two shorter edges which means that the equation has to be coupled with homogeneous Navier boundary conditions on these two edges and free boundary conditions on the other two edges. We refer to these boundary conditions as $(BC)$.

In Theorem \ref{t:Lax-Milgram} we show existence and uniqueness of weak solutions for the equation $\mathcal A u=f$ coupled with $(BC)$; in the same statement we also prove a regularity result based on elliptic regularity estimates by \cite{adn}.

The last part of the paper is devoted to the evolution problem for the plate which can be obtained with classical variational methods starting from the Lagrangian $\mathcal L$ given by a combination of the kinetic energy and the bending energy. Denoting by $\mathbb E_B(u)$ the bending energy corresponding to the configuration determined by the displacement function $u$ and by $M$ the linear mass density of the plate ($M$ is given by the total mass of the plate divided by its length $L$), we may write
\begin{equation*}
    \mathcal L(u):=\frac 12 \int_\Omega \frac{M}{2\ell} \left(\frac{\partial u}{\partial t}\right)^2 \, dxdy-\mathbb E_B(u) \, .
\end{equation*}

Then from $\mathcal L$ we obtain the following equation of motion
\begin{equation*}
  \frac{M}{2\ell} \frac{\partial^2 u}{\partial t^2}+\mathcal Au=0
\end{equation*}
coupled with the boundary conditions $(BC)$.

Special solutions of the equation of motion are the so-called stationary wave solutions. We recall how these solutions can be constructed starting from the eigenfunctions of $\mathcal A$ coupled with $(BC)$. Classical spectral theory implies that
this eigenvalue problem admits a discrete spectrum whose eigenvalues may be ordered in an increasing sequence diverging to $+\infty$. If $\lambda$ is an eigenvalue of $\mathcal A$ with corresponding eigenfunction $U_\lambda=U_\lambda(x,y)$, a stationary wave solution admits the following representation
$$
   u(x,y,t)=\sin(\omega_\lambda t) U_\lambda(x,y)
$$
where $\omega_\lambda$ is the ``angular velocity'' whose value is uniquely determined by $\lambda$ as we can see from Section \ref{s:vibrations}; then the frequencies of free vibration $\nu_\lambda$ are immediately obtained by writing $\nu_\lambda=\omega_\lambda/(2\pi)$.

The paper is organized as follows: Section \ref{s:orthotropic} is devoted to a general discussion on the basic notions about orthotropic materials, Section \ref{s:orthotropic-plate} is devoted to the construction of the model of orthotropic plate and to existence, uniqueness and regularity of solutions of the equilibrium equation under the action of a vertical load, and  Section \ref{s:vibrations} is devoted to the evolution problem, the eigenvalue problem, the construction of the stationary wave solutions and the corresponding frequencies of free vibration. In the final part of Section \ref{s:vibrations} we also suggest a possible model for a complete suspension bridge consisting of a system of three coupled equations, one for the deck and one for each of the two cables.


\section{Orthotropic materials} \label{s:orthotropic}

\subsection{Basic notions on the theory of anisotropic materials}
We recall some notations from the theory of linear elasticity. We denote by $\SI=(\sigma_{ij})$ the stress tensor and by ${\bf e}=(e_{ij})$  the strain tensor where $i,j\in \{1,2,3\}$ for both tensors.
We clarify that by strain tensor ${\bf e}$ we actually mean here the ``linearized strain tensor'', i.e.
\begin{equation} \label{eq:strain-tensor}
e_{ij}=\frac 12 \left(\frac{\partial u_i}{\partial
x_j}+\frac{\partial u_j}{\partial x_i}\right) \, , \qquad i,j\in
\{1,\dots,3\} \, .
\end{equation}
Here ${\bf u}=(u_1,u_2,u_3)$ is the vector field which defines the displacement at every point of the elastic body.

We recall that both ${\bm \sigma}$ and ${\bf e}$ are symmetric tensors, i.e. $\sigma_{ij}=\sigma_{ji}$ and $e_{ij}=e_{ji}$ for any $i,j\in \{1,2,3\}$.
See the book by Landau \& Lifshitz \cite{LaLi} for more details on these basic notions.

In the theory of linear elasticity, the tension of a material as a consequence of a deformation is proportional to the deformation itself and conversely the deformation of a material is proportional to the forces acting on it. This notion can be resumed by the generalized Hooke's law that states the existence of $81$ coefficients $(C_{ijkl})$ with $i,j,k,l\in \{1,2,3\}$ such that
\begin{equation} \label{eq:Hook-0}
  \sigma_{ij}=\sum_{k,l=1}^{3} C_{ijkl}\,  e_{kl}   \qquad \text{for any } i,j\in \{1,2,3\} \, .
\end{equation}
The elastic energy per unit of volume, as a function of the components of the strain tensor, can be implicitly characterized by
\begin{equation} \label{eq:implcit}
\frac{\partial \mathcal E}{\partial e_{ij}}=\sigma_{ij} \, , \ \ \text{for any } i,j\in \{1,2,,3\} \, , \quad   \mathcal E=0
  \ \ \ \text{when the material is undeformed} \, ,
\end{equation}
see \cite[Chapter 1, Paragraph 2]{LaLi} for more explanations on this question.

Then, by \eqref{eq:Hook-0}, we have that
\begin{equation} \label{eq:second-der}
  \frac{\partial^2 \mathcal E}{\partial e_{ij}\partial e_{kl}}=C_{ijkl} \qquad \text{for any } i,j,k,l\in \{1,2,3\}
\end{equation}
thus showing the following symmetry property
\begin{equation} \label{eq:symmetry-1}
  C_{ijkl}=C_{klij}  \qquad \text{for any } i,j,k,l\in \{1,2,3\}  \, .
\end{equation}
Combining \eqref{eq:implcit}, \eqref{eq:second-der} and \eqref{eq:symmetry-1} we obtain the explicit representation of the elastic energy per unit of volume:
\begin{equation} \label{eq:explicit}
 \mathcal E=\frac 12 \sum_{i,j,k,l=1}^{3} C_{ijkl} \, e_{ij} e_{kl} \, .
\end{equation}

Since the tensors $\SI$ and $\e$ are symmetric, actually each of them is completely determined by only six components and the related Hooke's law admits the following matrix representation
{\small
\begin{equation*}
  \begin{pmatrix}
    \sigma_{11} \\
    \sigma_{22} \\
    \sigma_{33} \\
    \sigma_{12} \\
    \sigma_{13} \\
    \sigma_{23}
  \end{pmatrix}
  =
  \begin{pmatrix}
    C_{1111} & C_{1122} & C_{1133} & C_{1112} & C_{1113} & C_{1123} \\
    C_{2211} & C_{2222} & C_{2233} & C_{2212} & C_{2213} & C_{2223} \\
    C_{3311} & C_{3322} & C_{3333} & C_{3312} & C_{3313} & C_{3323} \\
    C_{1211} & C_{1222} & C_{1233} & C_{1212} & C_{1213} & C_{1223} \\
    C_{1311} & C_{1322} & C_{1333} & C_{1312} & C_{1313} & C_{1323} \\
    C_{2311} & C_{2322} & C_{2333} & C_{2312} & C_{2313} & C_{2323}
  \end{pmatrix}
  \begin{pmatrix}
    e_{11} \\
    e_{22} \\
    e_{33} \\
    e_{12} \\
    e_{13} \\
    e_{23}
  \end{pmatrix} \, .
\end{equation*}
}
We denote by $C$ the $6\times 6$ matrix appearing above and we refer to it as the stiffness matrix.

\subsection{Orthotropic materials and their stiffening matrix} \label{a:1}
A material is said to be orthotropic if the corresponding stiffness matrix remains invariant under reflections with respect to three mutually orthogonal planes. In other words there exists three
orthogonal symmetry planes for which the stiffness matrix remains invariant under the corresponding reflections.

We collect in this section a number of well known facts coming from the theory of orthotropic materials; such facts are presented here in details for the reader convenience.

Let us consider two orthonormal coordinate systems $x_1x_2x_3$ and $x_1'x_2'x_3'$, one related to the other by transformations of the form ${\bf x'}=A {\bf x}$ including rotations and reflections where $A=(A_{ij})$ is a $3\times 3$ orthogonal matrix and where we put ${\bf x}=(x_1,x_2,x_3)$ and ${\bf x'}=(x_1',x_2',x_3')$.

If we have a symmetric tensor which is represented by a matrix $X$ in the coordinate system $x_1x_2x_3$ and we want to determine the corresponding matrix $X'$ in the new coordinate system $x_1'x_2'x_3'$, then one finds that $X'=AXA^T$.

Let us introduce the linear map $\mathcal L_A:{\rm Sym_3}\to {\rm Sym_3}$ defined by $\mathcal L_A(X)=AXA^T$ for any $X\in {\rm Sym_3}$, where ${\rm Sym_3}$ is the vector space of all $3\times 3$ symmetric matrices. As a base of ${\rm Sym_3}$ we choose the set of six matrices $\{X_1,\dots,X_6\}$ defined by
{\small
\begin{align} \label{eq:base}
\left\{
\begin{pmatrix}
1 & 0 & 0 \\
0 & 0 & 0 \\
0 & 0 & 0
\end{pmatrix} \, ,
\
\begin{pmatrix}
0 & 0 & 0 \\
0 & 1 & 0 \\
0 & 0 & 0
\end{pmatrix} \, ,
\
\begin{pmatrix}
0 & 0 & 0 \\
0 & 0 & 0 \\
0 & 0 & 1
\end{pmatrix} \, ,
\
\begin{pmatrix}
0 & 1 & 0 \\
1 & 0 & 0 \\
0 & 0 & 0
\end{pmatrix} \, ,
\
\begin{pmatrix}
0 & 0 & 1 \\
0 & 0 & 0 \\
1 & 0 & 0
\end{pmatrix} \, ,
\
\begin{pmatrix}
0 & 0 & 0 \\
0 & 0 & 1 \\
0 & 1 & 0
\end{pmatrix}
\right\} \, .
\end{align}
}
In this way, any matrix $X\in {\rm Sym_3}$ may be represented as $X=\sum_{i=1}^{6} \alpha_i X_i$. In the base $\{X_1,\dots,X_6\}$ the linear operator $\mathcal L_A$ may be represented by mean of a $6\times 6$ matrix in the following way
{\small
\begin{equation} \label{eq:big-matrix}
\begin{pmatrix}
  \alpha_1' \\
  \alpha_2' \\
  \alpha_3' \\
  \alpha_4' \\
  \alpha_5' \\
  \alpha_6'
\end{pmatrix}
=
\begin{pmatrix}
  A_{11}^2     & A_{12}^2     & A_{13}^2     & 2A_{11}A_{12}             & 2A_{11}A_{13}             & 2A_{12}A_{13} \\
  A_{21}^2     & A_{22}^2     & A_{23}^2     & 2A_{21}A_{22}             & 2A_{21}A_{23}             & 2A_{22}A_{23} \\
  A_{31}^2     & A_{32}^2     & A_{33}^2     & 2A_{31}A_{32}             & 2A_{31}A_{33}             & 2A_{32}A_{33} \\
  A_{11}A_{21} & A_{12}A_{22} & A_{13}A_{23} & A_{11}A_{22}+A_{12}A_{21} & A_{11}A_{23}+A_{13}A_{21} & A_{12}A_{23}+A_{13}A_{22} \\
  A_{11}A_{31} & A_{12}A_{32} & A_{13}A_{33} & A_{11}A_{32}+A_{12}A_{31} & A_{11}A_{33}+A_{13}A_{31} & A_{12}A_{33}+A_{13}A_{32} \\
  A_{21}A_{31} & A_{22}A_{32} & A_{23}A_{33} & A_{21}A_{32}+A_{22}A_{31} & A_{21}A_{33}+A_{23}A_{31} & A_{22}A_{33}+A_{23}A_{32}
\end{pmatrix}
\begin{pmatrix}
  \alpha_1 \\
  \alpha_2 \\
  \alpha_3 \\
  \alpha_4 \\
  \alpha_5 \\
  \alpha_6
\end{pmatrix}
\end{equation}
}
where we put $\mathcal L_A(X)=\sum_{i=1}^{6} \alpha_i' X_i$.

We denote by ${\mathbb A}$ the $6\times 6$ matrix defined in \eqref{eq:big-matrix} corresponding to the matrix transformation $A$.

Suppose now that ${\bf e}$ are ${\bf \sigma}$ are the strain and stress tensors of some linear elastic material. Once we fix an orthonormal coordinate system $x_1x_2x_3$, they can both be represented by two $3\times 3$ matrices. We recall that ${\bf e}$ are ${\bf \sigma}$ are symmetric tensors so that they can both be represented as linear combinations of the matrices $X_1,\dots,X_6$.


If we perform a coordinate transformation from a system $x_1x_2x_3$ to a new orthonormal system $x_1'x_2'x_3'$ with  corresponding transformation $A$, then denoting by
\begin{align*}
& \widetilde{\bm\sigma}=
\begin{pmatrix}
\sigma_{11} & \sigma_{22} & \sigma_{33} & \sigma_{12} & \sigma_{13} & \sigma_{23}
\end{pmatrix}^T \, ,
\qquad
\widetilde{\bm e}=
\begin{pmatrix}
e_{11} & e_{22} & e_{33} & e_{12} & e_{13} & e_{23}
\end{pmatrix}^T \, ,
\end{align*}
the vector columns of components of the stress and strain tensors in the $x_1x_2x_3$ coordinate system
and by
\begin{align*}
& \widetilde{\bm\sigma}'=
\begin{pmatrix}
\sigma_{11}' & \sigma_{22}' & \sigma_{33}' & \sigma_{12}' & \sigma_{13}' & \sigma_{23}'
\end{pmatrix}^T \, ,
\qquad
\widetilde{\bm e}'=
\begin{pmatrix}
e_{11}' & e_{22}' & e_{33}' & e_{12}' & e_{13}' & e_{23}'
\end{pmatrix}^T \, ,
\end{align*}
the vector columns of components of the stress and strain tensors in the $x_1'x_2'x_3'$ coordinate system,
then by \eqref{eq:big-matrix} we obtain
\begin{equation*}
  \widetilde{\bm\sigma}'=\mathbb A \widetilde{\bm\sigma} \qquad \text{and} \qquad
  \widetilde{\bm e}'=\mathbb A \widetilde{\bm e} \, .
\end{equation*}
Therefore, denoting now by $C$ and $C'$ the stiffness matrices of the material in the two coordinate systems, we infer
\begin{equation*}
 C'=\mathbb A C \mathbb A^{-1} \, .
\end{equation*}

Recalling the definition of orthotropy given at the beginning of this section, we introduce the reflections with respect to the three coordinate planes whose matrices are given by
{\small
\begin{equation*}
  A_1=
  \begin{pmatrix}
    -1 & 0 & 0 \\
     0 & 1 & 0 \\
     0 & 0 & 1
  \end{pmatrix} \, ,
  \qquad
  A_2=
  \begin{pmatrix}
     1 &  0 & 0 \\
     0 & -1 & 0 \\
     0 &  0 & 1
  \end{pmatrix} \, ,
  \qquad
  A_3=
  \begin{pmatrix}
     1 & 0 &  0 \\
     0 & 1 &  0 \\
     0 & 0 &
     -1
  \end{pmatrix} \, .
\end{equation*}
}
We denote by $\mathbb A_1, \mathbb A_2, \mathbb A_3$ the $6\times 6$ matrices corresponding to $A_1, A_2, A_3$ respectively, according with \eqref{eq:big-matrix}.
We observe that $\mathbb A_1, \mathbb A_2, \mathbb A_3$ are three diagonal matrices with a diagonal containing four times
the number $1$ and twice the number $-1$. We omit for simplicity the explicit representation of $\mathbb A_1, \mathbb A_2, \mathbb A_3$ and we proceed directly to the computation of the following matrices

{\small
\begin{equation*}
\mathbb A_1 C \mathbb A_1^{-1}
=
\begin{pmatrix}
    C_{1111} &  C_{1122} &  C_{1133} & -C_{1112} & -C_{1113} &  C_{1123} \\
    C_{2211} &  C_{2222} &  C_{2233} & -C_{2212} & -C_{2213} &  C_{2223} \\
    C_{3311} &  C_{3322} &  C_{3333} & -C_{3312} & -C_{3313} &  C_{3323} \\
   -C_{1211} & -C_{1222} & -C_{1233} &  C_{1212} &  C_{1213} & -C_{1223} \\
   -C_{1311} & -C_{1322} & -C_{1333} &  C_{1312} &  C_{1313} & -C_{1323} \\
    C_{2311} &  C_{2322} &  C_{2333} & -C_{2312} & -C_{2313} &  C_{2323}
  \end{pmatrix} \, ,
\end{equation*}
}
{\small
\begin{equation*}
\mathbb A_2 C \mathbb A_2^{-1}
=
\begin{pmatrix}
    C_{1111} &  C_{1122} &  C_{1133} & -C_{1112} &  C_{1113} & -C_{1123} \\
    C_{2211} &  C_{2222} &  C_{2233} & -C_{2212} &  C_{2213} & -C_{2223} \\
    C_{3311} &  C_{3322} &  C_{3333} & -C_{3312} &  C_{3313} & -C_{3323} \\
   -C_{1211} & -C_{1222} & -C_{1233} &  C_{1212} & -C_{1213} &  C_{1223} \\
    C_{1311} &  C_{1322} &  C_{1333} & -C_{1312} &  C_{1313} & -C_{1323} \\
   -C_{2311} & -C_{2322} & -C_{2333} &  C_{2312} & -C_{2313} &  C_{2323}
  \end{pmatrix} \, ,
\end{equation*}
}
{\small
\begin{equation*}
\mathbb A_3 C \mathbb A_3^{-1}
=
\begin{pmatrix}
    C_{1111} &  C_{1122} &  C_{1133} &  C_{1112} & -C_{1113} & -C_{1123} \\
    C_{2211} &  C_{2222} &  C_{2233} &  C_{2212} & -C_{2213} & -C_{2223} \\
    C_{3311} &  C_{3322} &  C_{3333} &  C_{3312} & -C_{3313} & -C_{3323} \\
    C_{1211} &  C_{1222} &  C_{1233} &  C_{1212} & -C_{1213} & -C_{1223} \\
   -C_{1311} & -C_{1322} & -C_{1333} & -C_{1312} &  C_{1313} &  C_{1323} \\
   -C_{2311} & -C_{2322} & -C_{2333} & -C_{2312} &  C_{2313} &  C_{2323}
  \end{pmatrix} \, .
\end{equation*}
}

The orthotropy condition implies that the matrix $C$ coincides simultaneously with all the three matrices $\mathbb A_1 C \mathbb A_1^{-1}$, $\mathbb A_2 C \mathbb A_2^{-1}$ and $\mathbb A_3 C \mathbb A_3^{-1}$ thus showing that the stiffness matrix $C$ is in the form
{\small
\begin{equation} \label{eq:stiff-ortho}
  \begin{pmatrix}
    C_{1111} & C_{1122} & C_{1133} & 0        & 0        & 0 \\
    C_{2211} & C_{2222} & C_{2233} & 0        & 0        & 0 \\
    C_{3311} & C_{3322} & C_{3333} & 0        & 0        & 0 \\
    0        & 0        & 0        & C_{1212} & 0        & 0 \\
    0        & 0        & 0        & 0        & C_{1313} & 0 \\
    0        & 0        & 0        & 0        & 0        & C_{2323}
  \end{pmatrix} \, .
\end{equation}
}

Another property of the stiffening matric $C$ is the following:

\begin{equation} \label{eq:C2323}
   C_{2323}=\frac{C_{2222}-C_{2233}-C_{3322}+C_{3333}}{2} \, .
\end{equation}

In order to prove \eqref{eq:C2323}, let us consider a transformation matrix $A$ in the form
{\small
\begin{equation*}
A=R_\theta:=
\begin{pmatrix}
1 & 0 & 0 \\
0 & \cos \theta & -\sin\theta \\
0 & \sin\theta & \cos\theta
\end{pmatrix} \, , \qquad \theta\in (-2\pi,2\pi) \, .
\end{equation*}
}
The corresponding matrix $\mathbb A$ introduced in \eqref{eq:big-matrix} becomes
\begin{equation} \label{eq:op-3}
\mathbb A=T_\theta:=\begin{pmatrix}
1 & 0                      & 0                     & 0           & 0           & 0              \\
0 & \cos^2 \theta          & \sin^2\theta          & 0           & 0           & -\sin(2\theta) \\
0 & \sin^2\theta           & \cos^2\theta          & 0           & 0           & \sin(2\theta)  \\
0 & 0                      & 0                     & \cos\theta  & -\sin\theta & 0              \\
0 & 0                      & 0                     & \sin\theta  & \cos\theta  & 0              \\
0 & \sin\theta \cos \theta & -\sin\theta\cos\theta & 0           & 0           & \cos(2\theta)
\end{pmatrix} \, .
\end{equation}
As explained above, the invariance with respect to the transformation $R_\theta$ implies that
\begin{equation} \label{eq:invariance-2323}
C=\mathbb AC\mathbb A^{-1}=T_\theta C T_\theta^{-1}
\end{equation}
where $C$ is the stiffening matrix. It can be easily deduced that the inverse matrix of $T_\theta$ is the matrix $T_{-\theta}$.

We are now interested to exploit the identity \eqref{eq:invariance-2323} relatively to the components placed at the sixth row and sixth column since this will produce the proof of \eqref{eq:C2323}.
Indeed, combining \eqref{eq:stiff-ortho},\eqref{eq:op-3} and \eqref{eq:invariance-2323} we obtain
\begin{equation*} 
C_{2323}=(\mathbb T_\theta C\mathbb T_\theta^{-1})_{66}=\frac{C_{2222}-C_{2233}-C_{3322}+C_{3333}}{2}\, \sin^2(2\theta)+C_{2323} \cos^2(2\theta)
\end{equation*}
and hence
\begin{equation*}
   C_{2323} [1-\cos^2(2\theta)]=\frac{C_{2222}-C_{2233}-C_{3322}+C_{3333}}{2}\, \sin^2(2\theta)  \, .
\end{equation*}
Choosing $\theta\in (0,2\pi)$, $\theta \neq \frac{\pi}{2}, \pi, \frac{3\pi}{2}$ and simplifying, this gives \eqref{eq:C2323}.

We observe that, differently from $C_{2323}$, the coefficients $C_{1212}$ and $C_{1313}$ are completely independent of the other ones as one can understand observing that proceeding again with the transformation $R_\theta$, the components of
$T_\theta C T_\theta^{-1}$ placed at the fourth/fifth row and fourth/fifth
column are given by $C_{1212}$ and $C_{1313}$ respectively, independently of the value assumed by the angle $\theta$, thus
producing no restrictions on $C_{1212}$ and $C_{1313}$.

\subsection{The elastic energy and the inverse Hooke's law} \label{ss:Ortho-Mat}

As shown in Section \ref{a:1}, the Hooke's law for an orthotropic material reads
{\small
\begin{equation} \label{eq:Hook-1}
\begin{pmatrix}
\sigma_{11} \\
\sigma_{22} \\
\sigma_{33} \\
\sigma_{12} \\
\sigma_{13} \\
\sigma_{23} \\
\end{pmatrix}
=
\begin{pmatrix}
C_{1111} & C_{1122} & C_{1133} & 0 & 0 & 0 \\
C_{1122} & C_{2222} & C_{2233} & 0 & 0 & 0 \\
C_{1133} & C_{2233} & C_{3333} & 0 & 0 & 0 \\
0 & 0 & 0 & C_{1212} & 0 & 0 \\
0 & 0 & 0 & 0 & C_{1313} & 0 \\
0 & 0 & 0 & 0 & 0 & C_{2323}
\end{pmatrix}
\begin{pmatrix}
e_{11} \\
e_{22} \\
e_{33} \\
e_{12} \\
e_{13} \\
e_{23} \\
\end{pmatrix} \, .
\end{equation}
}

In order to apply \eqref{eq:explicit}, we collect the coefficients $C_{ijlm}$ introduced in \eqref{eq:Hook-0} in a $9\times 9$ matrix that we denote by ${\bf C}_{9\times 9}$. Invoking \eqref{eq:Hook-1}, we see that in the special case of an orthotropic material, the matrix ${\bf C}_{9\times 9}$ admits the following representation in terms of $3\times 3$ blocks
\begin{equation} \label{eq:blocks}
  \left(
  \begin{tabular}{c|c|c}
    ${\bf B_1}$ & $\bf 0$ & $\bf 0$ \\
    \hline
    $\bf 0$ & ${\bf B_2}$ & $\bf 0$ \\
    \hline
    $\bf 0$ & $\bf 0$ & ${\bf B_3}$ \\
  \end{tabular}
  \right)
\end{equation}
where
{\small
\begin{equation*}
  {\bf B_1}=\begin{pmatrix}
C_{1111} & C_{1122} & C_{1133} \\
C_{1122} & C_{2222} & C_{2233} \\
C_{1133} & C_{2233} & C_{3333}
\end{pmatrix} \, ,
\quad
 {\bf B_2}=\begin{pmatrix}
C_{1212} & 0        & 0 \\
0        & C_{1313} & 0 \\
0        & 0        & C_{2323}
\end{pmatrix} \, ,
\quad
 {\bf B_3}=\begin{pmatrix}
C_{2121} & 0        & 0 \\
0        & C_{3131} & 0 \\
0        & 0        & C_{3232}
\end{pmatrix}
\end{equation*}
and $\bf 0$ denotes the $3\times 3$ null matrix. From the symmetry of $\SI$ and $\e$ we deduce that
\begin{equation} \label{eq:symmetry-2}
  C_{1212}=C_{2121} \, , \qquad C_{1313}=C_{3131} \, , \qquad  C_{2323}=C_{3232}  \, ,
\end{equation}
}
and in particular ${\bf B_2}={\bf B_3}$.

If we replace the usual representations of $\SI$ and $\e$ as $3\times 3$ matrices with the following ones as vector columns of $9$ components
\begin{align*} 
& \SI=
\begin{pmatrix}
\sigma_{11} & \sigma_{22} & \sigma_{33} & \sigma_{12} &
\sigma_{13} & \sigma_{23} & \sigma_{21} & \sigma_{31} &
\sigma_{32}
\end{pmatrix}^T \, ,  \\
& \e =
\begin{pmatrix}
e_{11} & e_{22} & e_{33} & e_{12} & e_{13} & e_{23} & e_{21} &
e_{31} & e_{32}
\end{pmatrix}
^T \, ,
\end{align*}
by \eqref{eq:explicit}, \eqref{eq:Hook-1} and \eqref{eq:blocks} we deduce that
\begin{equation*} 
\SI={\bf C}_{9\times 9} \, \e  \qquad \text{and} \qquad \mathcal E=\frac 12 \, \e^T {\bf C}_{9\times 9} \, \e \, .
\end{equation*}
In particular by \eqref{eq:blocks}, \eqref{eq:symmetry-2} and the symmetry of $\e$, we have
\begin{align}\label{eq:elastic-energy}
  \mathcal E & =\frac 12 \,
  \e_{{\rm diag}}^T  \,
  {\bf B_1} \, \e_{{\rm diag}}
  +\frac{C_{1212}\, e_{12}^2
  +C_{1313}\, e_{13}^2
  +C_{2323}\, e_{23}^2}2
  +\frac{C_{2121}\, e_{21}^2
  +C_{3131}\, e_{31}^2
  +C_{3232}\, e_{32}^2}2 \\[10pt]
  & \notag  = \frac 12 \,
  \e_{{\rm diag}}^T  \,
  {\bf B_1} \, \e_{{\rm diag}}
  +C_{1212}\, e_{12}^2
  +C_{1313}\, e_{13}^2
  +C_{2323}\, e_{23}^2
\end{align}
with $\e_{{\rm diag}}=\begin{pmatrix} e_{11} & e_{22} & e_{33} \end{pmatrix}^T$.

Consider now the inverse of identity \eqref{eq:Hook-1} in such a way that the strain tensor is expressed in terms of the stress tensor,
\begin{equation} \label{eq:Hook-2}
\begin{pmatrix}
e_{11} \\
e_{22} \\
e_{33} \\
e_{12} \\
e_{13} \\
e_{23} \\
\end{pmatrix}
=
\begin{pmatrix}
\frac 1{E_{1}} & -\frac{\nu_{21}}{E_{2}} & -\frac{\nu_{31}}{E_{3}} & 0 & 0 & 0 \\
-\frac{\nu_{12}}{E_{1}} & \frac 1{E_{2}} & -\frac{\nu_{32}}{E_{3}} & 0 & 0 & 0 \\
-\frac{\nu_{13}}{E_{1}} & -\frac{\nu_{23}}{E_{2}} & \frac 1{E_{3}} & 0 & 0 & 0 \\
0 & 0 & 0 & \frac 1{2\mu_{12}} & 0 & 0 \\
0 & 0 & 0 & 0 & \frac 1{2\mu_{13}} & 0 \\
0 & 0 & 0 & 0 & 0 & \frac 1{2\mu_{23}}
\end{pmatrix}
\begin{pmatrix}
\sigma_{11} \\
\sigma_{22} \\
\sigma_{33} \\
\sigma_{12} \\
\sigma_{13} \\
\sigma_{23} \\
\end{pmatrix}
\end{equation}
where the constants $E_{1}, E_{2}, E_{3}$ are known as Young
moduli, one for each direction and the constants
$\nu_{ij}$, $i,j\in \{1,2,3\}$, $i\neq j$ are known as Poisson
ratios. Finally $\mu_{12}, \mu_{13}, \mu_{23}$ are known as shear moduli
or moduli of rigidity and their indices coincide with the corresponding components of the stress and strain tensors obviously involved by these coefficients:
\begin{equation*}
  e_{12}=(2\mu_{12})^{-1} \, \sigma_{12} \, , \qquad e_{13}=(2\mu_{13})^{-1} \, \sigma_{13} \, ,\qquad e_{23}=(2\mu_{23})^{-1} \, \sigma_{23} \, .
\end{equation*}

Let us denote by $S$ the $6\times 6$ matrix appearing in \eqref{eq:Hook-2}. Comparing \eqref{eq:Hook-1} with \eqref{eq:Hook-2}, it appears clear that $S=C^{-1}$.
The symmetry properties of the stiffening matrix $C$, see
\eqref{eq:Hook-1}, implies symmetry of $S$ so that we may write
\begin{equation} \label{eq:id0}
\frac{\nu_{21}}{E_{2}}=\frac{\nu_{12}}{E_{1}} \, , \qquad
\frac{\nu_{31}}{E_{3}}=\frac{\nu_{13}}{E_{1}} \, , \qquad
\frac{\nu_{32}}{E_{3}}=\frac{\nu_{23}}{E_{2}} \, .
\end{equation}
We observe that by \eqref{eq:Hook-2} we easily deduce that in the
case of a one-dimensional tension state parallel to the $x_1$-axis, i.e.
the only component of $\SI$ different from zero is $\sigma_{11}$,
we have that
\begin{equation*}
\nu_{12}=-\tfrac{e_{22}}{e_{11}} \quad \text{and} \quad
\nu_{13}=-\tfrac{e_{33}}{e_{11}} \, .
\end{equation*}
Similarly, choosing first a one-dimensional tension state parallel
to the $x_2$-axis and then a one-dimensional tension state parallel to the $x_3$-axis, we infer
\begin{equation*}
\nu_{21}=-\tfrac{e_{11}}{e_{22}} \, , \qquad \nu_{23}=-\tfrac{e_{33}}{e_{22}} \, , \qquad
\nu_{31}=-\tfrac{e_{11}}{e_{33}} \, , \qquad \nu_{32}=-\tfrac{e_{22}}{e_{33}} \, .
\end{equation*}
In other words, if we consider the Poisson ratio $\nu_{ij}$, the index $i$ (the first one) represents the direction of the one-dimensional stress and the index $j$ (the second one) represents the direction of the transversal deformation.
This explanation clarifies the meaning of the Poisson ratios and the notation used in \eqref{eq:Hook-2}.

We want now to represent the components of the matrix $C$ in terms of the coefficients $E_i$, $\nu_{ij}$, $\mu_{ij}$
introduced in \eqref{eq:Hook-2}, by proceeding with the inversion of $S$ and then by comparing the components of its inverse with the coefficients $C_{ijkl}$.

Let us put
\begin{equation} \label{eq:def-delta}
\delta:={\rm det}
\begin{pmatrix}
  \frac{1}{E_1}           & -\frac{\nu_{21}}{E_2} & -\frac{\nu_{31}}{E_3} \\
  -\frac{\nu_{12}}{E_1}   & \frac{1}{E_2}         & -\frac{\nu_{32}}{E_3} \\
  -\frac{\nu_{13}}{E_1}   & -\frac{\nu_{23}}{E_2} & \frac{1}{E_3}
\end{pmatrix}
=\frac{1-\nu_{12}\nu_{21}-\nu_{13}\nu_{31}-\nu_{23}\nu_{32}-2\nu_{12}\nu_{23}\nu_{31}}{E_1
E_2 E_3}
\end{equation}
where we exploited \eqref{eq:id0} to show that $\nu_{13}\nu_{21}\nu_{32}=\nu_{12}\nu_{23}\nu_{31}$.
In this way we may write $\det(S)=\frac{\delta}{8\mu_{12}\mu_{13}\mu_{23}}$.

Using this notation we obtain
\begin{equation} \label{eq:S-1}
C=S^{-1}=
\begin{pmatrix}
\frac{1-\nu_{23}\nu_{32}}{\delta \, E_2E_3} & \frac{\nu_{12}+\nu_{13}\nu_{32}}{\delta \, E_1 E_3} &
\frac{\nu_{13}+\nu_{12}\nu_{23}}{\delta \, E_1 E_2}
& 0 & 0 & 0 \\
\frac{\nu_{21}+\nu_{31}\nu_{23}}{\delta \, E_2 E_3} & \frac{1-\nu_{13}\nu_{31}}{\delta \, E_1 E_3} &
\frac{\nu_{23}+\nu_{13}\nu_{21}}{\delta \, E_1 E_2}
& 0 & 0 & 0 \\
\frac{\nu_{31}+\nu_{21}\nu_{32}}{\delta \, E_2 E_3} & \frac{\nu_{32}+\nu_{31}\nu_{12}}{\delta \, E_1 E_3}
& \frac{1-\nu_{12}\nu_{21}}{\delta \, E_1 E_2}
& 0 & 0 & 0 \\
0 & 0 & 0 & 2\mu_{12} & 0 & 0 \\
0 & 0 & 0 & 0 & 2\mu_{13} & 0 \\
0 & 0 & 0 & 0 & 0 & 2\mu_{23}
\end{pmatrix} \, .
\end{equation}

\subsection{Orthotropic materials with a one-dimensional reinforcement} \label{ss:p-iso}
Let us consider an orthotropic material which a one-dimensional reinforcement in the $x_1$ direction and an isotropic behavior in the $x_2$, $x_3$ variables. If we look at \eqref{eq:Hook-2}, this means that
\begin{equation} \label{eq:id1}
E_2=E_3 \, , \quad \nu_{21}=\nu_{31} \, , \quad \nu_{12}=\nu_{13}
\, , \quad \mu_{12}=\mu_{13} \, , \quad \nu_{23}=\nu_{32} \, .
\end{equation}

Recalling \eqref{eq:C2323} and exploiting \eqref{eq:id1}, we see that the elastic properties of the material are uniquely determined by only five constants:
\begin{equation} \label{eq:five-const}
E_1, \ E_2, \  \nu_{12}, \ \nu_{23}, \ \mu_{12} \, .
\end{equation}
Indeed, \eqref{eq:C2323} and \eqref{eq:id1} imply $C_{2222}=C_{3333}$, $C_{2233}=C_{3322}$ so that
\begin{equation} \label{eq:conjecture}
\mu_{23}=\tfrac{C_{2323}}2=\tfrac 12(C_{2222}-C_{2233}) \, .
\end{equation}
Now, exploiting the identities in \eqref{eq:id0}, \eqref{eq:id1}, \eqref{eq:conjecture}, we can write the matrix $C$ in the form
{\small
\begin{equation} \label{eq:rigidity-matrix}
C=
\begin{pmatrix}
\frac{1-\nu_{23}^2}{\delta \, E_2^2} & \frac{\nu_{12}(1+\nu_{23})}{\delta \, E_1 E_2} & \frac{\nu_{12}(1+\nu_{23})}{\delta \, E_1 E_2}
& 0 & 0 & 0 \\[7pt]
\frac{\nu_{12}(1+\nu_{23})}{\delta \, E_1 E_2} & \frac{E_1-E_2 \, \nu_{12}^2}{\delta \, E_1^2 E_2} & \frac{E_1 \nu_{23}+E_2 \, \nu_{12}^2}{\delta \, E_1^2 E_2}
& 0 & 0 & 0 \\[7pt]
\frac{\nu_{12}(1+\nu_{23})}{\delta \, E_1 E_2} & \frac{E_1 \, \nu_{23}+E_2 \, \nu_{12}^2}{\delta \, E_1^2 E_2} &  \frac{E_1-E_2 \, \nu_{12}^2}{\delta \, E_1^2 E_2}
& 0 & 0 & 0 \\[7pt]
0 & 0 & 0 & 2\mu_{12} & 0 & 0 \\[7pt]
0 & 0 & 0 & 0 & 2\mu_{12} & 0 \\[7pt]
0 & 0 & 0 & 0 & 0 & \frac{E_1(1-\nu_{23})-2E_2 \, \nu_{12}^2}{\delta \, E_1^2 E_2}
\end{pmatrix} \, .
\end{equation}
}

\section{The model of a plate with a one-dimensional reinforcement} \label{s:orthotropic-plate}
We present in this section the model of a plate of length $L$, width $2\ell$ and
thickness $d$. We choose a coordinate system in such a way that the plate is described by
the set $(0,L)\times (-\ell,\ell)\times \left(-\frac
d2,\frac d2 \right)$. We use for this coordinates the usual $x,y,z$ notation in place of the $x_1, x_2, x_3$ notation
used in Section \ref{s:orthotropic}.

We assume that the plate is made of an
orthotropic material with a one-dimensional reinforcement in the
$x$ direction. We assume the validity of the classical
constitutive assumptions for the displacement of a plate, see
\cite[Paragraph 11]{LaLi}:

\begin{itemize}
\item the displacement of the midway surface is only vertical and it is described by a function $u=u(x,y)$ with $(x,y)\in (0,L)\times (-\ell,\ell)$;

\item the third component of the displacement vector ${\bf u}=(u_1,u_2,u_3)$ only depends on $x$ and $y$ and with sufficient accuracy we may assume that $u_3(x,y)=u(x,y)$ for any $(x,y)\in (0,L)\times (-\ell,\ell)$;

\item the components $\sigma_{13}, \sigma_{23}, \sigma_{33}$ of the stress tensor vanish everywhere in the plate.
\end{itemize}

We now compute the elastic energy per unit of volume in a configuration corresponding to a
generic vertical displacement $u$ of the midway surface. By \eqref{eq:strain-tensor} and the above constitutive conditions we obtain
\begin{equation*}
u_{1}=-z \frac{\partial u}{\partial x} \, , \qquad u_{2}=-z \frac{\partial u}{\partial y} \, , \qquad u_3=u \, ,
\end{equation*}
and, in turn,
\begin{equation} \label{eq:e11-e22}
e_{11}=-z \frac{\partial^2 u}{\partial x^2} \, , \qquad e_{22}=-z \frac{\partial^2 u}{\partial y^2} \, , \qquad
e_{12}=-z\frac{\partial^2 u}{\partial x\partial y} \, , \qquad e_{13}=0 \, , \qquad e_{23}=0 \, ,
\end{equation}
see \cite[Paragraph 11]{LaLi} for more details. Finally, condition $\sigma_{33}=0$ combined with \eqref{eq:Hook-1}, \eqref{eq:id0} and \eqref{eq:rigidity-matrix}, yields
\begin{align} \label{eq:e33}
e_{33}&=-\frac{\nu_{12}E_1 (1+\nu_{23})e_{11}+(E_1 \nu_{23}+E_2 \nu_{12}^2)e_{22}}{E_1-E_2 \nu_{12}^2} \\[7pt]
\notag & =\frac{\nu_{12}E_1 (1+\nu_{23})}{E_1-E_2 \nu_{12}^2} \, z \frac{\partial^2 u}{\partial x^2}+\frac{E_1 \nu_{23}+E_2 \nu_{12}^2}{E_1-E_2 \nu_{12}^2}
\, z \frac{\partial^2 u}{\partial y^2} \, .
\end{align}
Replacing \eqref{eq:e11-e22} and \eqref{eq:e33} into \eqref{eq:elastic-energy} we obtain
\begin{equation*}
\mathcal E=\frac 12 \left(K_{11} e_{11}^2+K_{22}
e_{22}^2+2K_{1122} \, e_{11}e_{22}+2K_{1212} e_{12}^2\right)
\end{equation*}
where
\begin{align*}
& K_{11}=\frac{(1+\nu_{23})[E_1(1-\nu_{23})-2E_2 \nu_{12}^2]}{\delta E_2^2 (E_1-E_2 \nu_{12}^2)} \, , \qquad
K_{22}=\frac{(1+\nu_{23})[E_1(1-\nu_{23})-2E_2 \nu_{12}^2]}{\delta E_1 E_2(E_1-E_2 \nu_{12}^2)}  \, , \\[10pt]
& K_{1122}=\frac{\nu_{12}(1+\nu_{23})[E_1(1-\nu_{23})-2E_2 \nu_{12}^2]}{\delta E_1 E_2 (E_1-E_2\nu_{12}^2)} \, , \qquad
K_{1212}=2\mu_{12} \, .
\end{align*}
Putting $\mathcal K=K_{22}$ we have that
\begin{align} \label{eq:mathcal-E}
\mathcal E & =\frac 12 \left(\frac{E_1 \mathcal K}{E_2} \,  e_{11}^2+\mathcal K e_{22}^2+2\nu_{12} \mathcal K \, e_{11}e_{22}+4\mu_{12} e_{12}^2\right) \\[7pt]
\notag & =\frac{z^2}2 \left[\frac{E_1 \mathcal K}{E_2}
\left(\frac{\partial^2 u}{\partial x^2}\right)^2+\mathcal
K\left(\frac{\partial^2 u}{\partial y^2}\right)^2 +2\nu_{12}
\mathcal K \frac{\partial^2 u}{\partial x^2}\frac{\partial^2
u}{\partial y^2}+4\mu_{12} \left(\frac{\partial^2 u}{\partial x
\partial y}\right)^2 \right] \, .
\end{align}
We observe that by \eqref{eq:id0}, \eqref{eq:def-delta} and some computations we may write $\mathcal K$ in a more elegant way:
\begin{equation} \label{eq:write-K}
  \mathcal K=\frac{E_2}{1-\nu_{12}\nu_{21}}  \, .
\end{equation}

Looking at \eqref{eq:mathcal-E} and \eqref{eq:write-K}, we see that the elastic coefficients that completely determines the bending energy of the plate corresponding to a generic displacement $u$ are $E_1$, $E_2$, $\nu_{12}$ and $\mu_{12}$ while no dependence on the Poisson ratio $\nu_{23}$ occurs.

The total bending energy of the deformed plate in term of the
vertical displacement assumes the form
\begin{equation} \label{eq:total-energy}
\mathbb E_B(u)=\frac{d^3}{24} \int_\Omega \left[\frac{E_1
\mathcal K}{E_2} \left(\frac{\partial^2 u}{\partial
x^2}\right)^2+\mathcal K\left(\frac{\partial^2 u}{\partial
y^2}\right)^2 +2\nu_{12} \mathcal K \frac{\partial^2 u}{\partial
x^2}\frac{\partial^2 u}{\partial y^2}+4\mu_{12}
\left(\frac{\partial^2 u}{\partial x
\partial y}\right)^2 \right] dxdy
\end{equation}
where we put $\Omega=(0,L)\times (-\ell,\ell)$.

As we pointed out at the end of Section \ref{ss:p-iso}, the
value of the coefficient $C_{1212}$, and hence of $\mu_{12}$,
is completely independent from the other coefficients so that the
choice of $\mu_{12}$ is free. It appears reasonable to our purposes to assume that
\begin{equation} \label{eq:ass-mu12}
\mu_{12}=\frac{\mathcal K (1-\nu_{12})}2 \, ,
\end{equation}
in complete accordance with the classical theory of plates in the isotropic setting in which $\mathcal K=\frac{E}{1-\nu^2}$
and $\mu=\frac{\mathcal K(1-\nu)}{2}=\frac{E}{2(1+\nu)}$, where $E$ and $\nu$ are the Young modulus and Poisson ratio of the isotropic material respectively and $\mu$ is one of the two Lam\'e coefficients usually known as modulus of rigidity of the material, see \cite[Chapter 1, Section 5, (5.9)]{LaLi}.

In this way, by \eqref{eq:total-energy}, we may write the bending
energy of the plate in the form
\begin{equation} \label{eq:total-energy-2}
\mathbb E_B(u)=\frac{d^3 \mathcal K}{24} \int_\Omega
\left[\frac{E_1}{E_2} \left(\frac{\partial^2 u}{\partial
x^2}\right)^2+\left(\frac{\partial^2 u}{\partial y^2}\right)^2
+2\nu_{12} \frac{\partial^2 u}{\partial x^2}\frac{\partial^2
u}{\partial y^2}+2(1-\nu_{12}) \left(\frac{\partial^2 u}{\partial
x
\partial y}\right)^2 \right] dxdy \, .
\end{equation}
Denoting by $D^2 u$ the Hessian matrix of $u$, introducing the notation
\begin{equation*}
  D^2 u:D^2 v=u_{xx} v_{xx}+2u_{xy}v_{xy}+u_{yy}v_{yy} \quad
\text{and} \quad |D^2 u|^2=u_{xx}^2+2u_{xy}^2+v_{yy}^2 \,
\end{equation*}
and defining
\begin{equation} \label{eq:write-k}
 \kappa=\frac{E_1-E_2}{E_2} \, ,
\end{equation}
 we may write
\begin{equation} \label{eq:E-B}
  \mathbb E_B(u)=\frac{d^3 \mathcal K}{24} \int_\Omega \left[\nu_{12} |\Delta u|^2+(1-\nu_{12})|D^2 u|^2+\kappa\,  u_{xx}^2   \right] dxdy \, .
\end{equation}
We observe that since the plate is reinforced in the $x$ direction, we assume that
\begin{equation} \label{eq:ipotesi-Young}
E_1>E_2
\end{equation}
which in turn implies $\kappa>0$.

In the remaining part of the paper, for the Poisson ratio $\nu_{12}$ we use the simpler notation
\begin{equation} \label{eq:write-nu}
    \nu=\nu_{12} \, .
\end{equation}
According with the theory of isotropic materials, we assume that
\begin{equation} \label{eq:ipotesi-Poisson}
0<\nu<\frac 12 \, .
\end{equation}

We now introduce a suitable functional space for the energy $\mathbb E_B$. As explained in the introduction, our main purpose is to describe the static and dynamic behavior of the deck of a bridge by mean of a plate model. For this reason, we may assume that the deck is hinged at the two vertical edges of the rectangle
$\Omega$ and free on the two horizontal edges of the same rectangle. Hence, a reasonable choice for the functional subspace of the Sobolev space $H^2(\Omega)$ is
\begin{equation} \label{eq:def-H2*}
H^2_*(\Omega):=\{w\in H^2(\Omega): w=0 \ \text{on} \ \{0,L\}\times
(-\ell,\ell)\} \, ,
\end{equation}
see \cite[Section 3]{FeGa}.

Thanks to the Intermediate Derivatives Theorem, see
\cite[Theorem 4.15]{Adams}, the space $H^2(\Omega)$ is a Hilbert
space if endowed with the scalar product
$$
(u,v)_{H^2}:=\int_\Omega\left(D^2u:D^2v+uv\right)\,
dxdy\qquad \text{for all } u,v\in H^2(\Omega)  \, .
$$

On the closed subspace $H^2_*(\Omega)$ it is possible to define an
alternative scalar product naturally related to the functional $\mathbb E_B$ as explained in the next proposition.

\begin{proposition} \label{l:equivalence}
Assume \eqref{eq:ass-mu12}, \eqref{eq:ipotesi-Young} and \eqref{eq:ipotesi-Poisson}. On
the space $H^2_*(\Omega)$ the two norms
$$
u\mapsto\|u\|_{H^2}\, ,\quad
u\mapsto\|u\|_{H^2_*}:=\left\{\int_\Omega \left[\nu |\Delta u|^2+(1-\nu)|D^2 u|^2+\kappa\,  u_{xx}^2   \right] dxdy\right\}^{1/2}
$$
are equivalent. Therefore, $H^2_*(\Omega)$ is a Hilbert space when
endowed with the scalar product
\begin{equation*} 
(u,v)_{H^2_*}:=\int_\Omega \left[ \nu \Delta u \Delta v+(1-\nu)D^2 u:D^2 v+\kappa \, u_{xx}v_{xx} \right] dxdy \, .
\end{equation*}
\end{proposition}

The proof can be obtained by proceeding as in the proof of
\cite[Lemma 4.1]{FeGa} or proceeding directly by combining the Poincar\'e inequality and the classical $H^2$-regularity estimate for the Laplacian.

Next, if we denote by $f$ an external vertical load per unit of
surface and if $u$ is the deflection of the plate in the vertical
direction, by \eqref{eq:total-energy-2} we
have that the total energy $\mathbb E_T$ of the plate becomes
\begin{equation} \label{eq:energy-total}
\mathbb E_T(u)=\frac{d^3 \mathcal K}{24} \int_\Omega
\left(\nu|\Delta u|^2+(1-\nu)|D^2 u|^2+\kappa\,
u_{xx}^2\right) dxdy-\int_\Omega fu \, dxdy \, .
\end{equation}
If $u\in C^4(\overline \Omega)\cap H^2_*(\Omega)$ and $v\in H^2_*(\Omega)$, by \cite[Proposition 5]{Chasman} and some
calculation, we infer
\begin{align*}
& \nu \int_\Omega \Delta u \Delta v \,
dxdy+(1-\nu)\int_\Omega D^2 u:D^2 v \, dxdy \\[6pt]
& = \int_0^L [\nu u_{xx}(x,\ell)+u_{yy}(x,\ell)]v_y(x,\ell)\,
dx-\int_0^L [\nu
u_{xx}(x,-\ell)+u_{yy}(x,-\ell)]v_y(x,-\ell)\, dx \\[6pt]
& -\!\!\int_0^L
[u_{yyy}(x,\ell)+(2-\nu)u_{xxy}(x,\ell)]v(x,\ell) \,
dx\!+\!\!\int_0^L
[u_{yyy}(x,-\ell)+(2-\nu)u_{xxy}(x,-\ell)]v(x,-\ell) \, dx
\\[6pt]
& \quad +\int_{-\ell}^\ell
[u_{xx}(L,y)v_x(L,y)-u_{xx}(0,y)v_x(0,y)]dy+\int_\Omega \Delta^2 u
\, v \, dxdy
\end{align*}
and
\begin{align*}
  \int_\Omega u_{xx}\, v_{xx} \, dxdy & = \int_{-\ell}^{\ell} u_{xx}(L,y)v_x(L,y)\, dy
  -\int_{-\ell}^{\ell} u_{xx}(0,y)v_x(0,y)\, dy+\int_\Omega u_{xxxx} \, v\, dxdy \, .
\end{align*}
Therefore, if $u\in C^4(\Omega)\cap H^2_*(\Omega)$ is a critical point of the functional $\mathbb E_T$ then it is a classical solution of the problem
\begin{equation} \label{eq:model-plate}
\begin{cases}
\frac{d^3\mathcal K}{12} \left(\Delta^2 u+\kappa\frac{\partial^4
u}{\partial x^4}\right)=f & \qquad
\text{in } \Omega \, , \\[6pt]
u(0,y)=u_{xx}(0,y)=u(L,y)=u_{xx}(L,y)=0 & \qquad \text{for }
y\in(-\ell,\ell) \, , \\[6pt]
u_{yy}(x,\pm \ell)+\nu u_{xx}(x,\pm \ell)=0 & \qquad \text{for }
x\in (0,L) \, , \\[6pt]
u_{yyy}(x,\pm \ell)+(2-\nu)u_{xxy}(x,\pm \ell)=0 & \qquad
\text{for } x\in (0,L) \, .
\end{cases}
\end{equation}
Problem \eqref{eq:model-plate} represents the model for a plate made of an orthotropic material with a one-dimensional reinforcement in the $x$ direction subject to vertical load $f$ per unit of surface.

\begin{remark} \label{r:1}
We observe that recalling \eqref{eq:id0} and \eqref{eq:write-K}, the fourth order equation in  \eqref{eq:model-plate} may be written in a different way, more familiar in the theory of orthotropic plates:
{\small
\begin{equation*}
 \frac{E_1 d^3}{12(1-\nu_{12} \nu_{21})} \, \frac{\partial^4 u}{\partial x^4}
 +\frac{E_2 d^3}{6(1-\nu_{12} \nu_{21})} \, \frac{\partial^4 u}{\partial x^2\partial y^2}
 +\frac{E_2 d^3}{12(1-\nu_{12} \nu_{21})} \, \frac{\partial^4 u}{\partial y^4}=f \, ,
\end{equation*}
}
see for example \cite[Chapter 2]{DesignManual}.
\end{remark}

Let us denote by $\mathcal H(\Omega)$ the dual space of $H^2_*(\Omega)$. We state here the following result taken from
\cite{Suspension}:

\begin{theorem} \label{t:Lax-Milgram}
Assume \eqref{eq:ass-mu12}, \eqref{eq:ipotesi-Young} and \eqref{eq:ipotesi-Poisson} and let $f\in \mathcal
H(\Omega)$. Then the following conclusions hold true:

\begin{itemize}

\item[$(i)$] there exists a unique $u\in H^2_*(\Omega)$ such
that
\begin{equation*}
\frac{d^3\mathcal K}{12} (u,v)_{H^2_*}=\ _{\mathcal
H(\Omega)}\langle f,v\rangle_{H^2_*(\Omega)} \qquad \text{for any
} v\in H^2_*(\Omega) \, ;
\end{equation*}

\item[$(ii)$] $u$ is the unique minimum point of the convex functional
$$
  \mathbb E_T(u)=\frac 12 (u,u)_{H^2_*}- \ _{\mathcal
H(\Omega)}\langle f,u\rangle_{H^2_*(\Omega)}\, ;
$$

\item[$(iii)$] if $f\in W^{k,p}(\Omega)$ for some $1<p<\infty$ and $k\in \N\cup\{0\}$, where we put $W^{0,p}(\Omega):=L^p(\Omega)$, then $u\in W^{k+4,p}(\Omega)$.

\end{itemize}
\end{theorem}

The proof of Theorem \ref{t:Lax-Milgram} is based on the Lax-Milgram Theorem and classical elliptic regularity estimates.
For the proof of this result we quote \cite{Suspension}.


\section{Eigenvalues and frequencies of free vibration} \label{s:vibrations}

This section is devoted to the evolution equation for the orthotropic plate with a particular attention for stationary wave solutions and for the related frequencies of vibration.

As anticipated in the introduction, the frequencies of vibration are closely related to the following eigenvalue problem
\begin{equation} \label{eq:eigenvalue-original}
\begin{cases}
\frac{d^3 \mathcal K}{12}\left(\Delta^2 u+\kappa \frac{\partial^4 u}{\partial x^4}\right)=\lambda u & \qquad
\text{in } \Omega=(0,L)\times (-\ell,\ell) \, , \\[6pt]
u(0,y)=u_{xx}(0,y)=u(L,y)=u_{xx}(L,y)=0 & \qquad \text{for }
y\in(-\ell,\ell) \, , \\[6pt]
u_{yy}(x,\pm \ell)+\nu u_{xx}(x,\pm \ell)=0 & \qquad \text{for }
x\in (0,L) \, , \\[6pt]
u_{yyy}(x,\pm \ell)+(2-\nu)u_{xxy}(x,\pm \ell)=0 &  \qquad
\text{for } x\in (0,L) \, .
\end{cases}
\end{equation}

Problem \eqref{eq:eigenvalue-original} admits a sequence of eigenvalues
\begin{equation} \label{eq:eig-plate}
  0<\lambda_1\le \lambda_2 \le \dots \le \lambda_m \le \dots
\end{equation}
diverging to $+\infty$. Indeed \eqref{eq:eigenvalue-original} admits the following weak formulation
\begin{equation*}
  \frac{ d^3 \mathcal K}{12} (u,v)_{H^2_*}=\lambda (u,v)_{L^2}  \qquad \text{for any } v\in H^2_*(\Omega)
\end{equation*}
so that compact embedding $H^2_*(\Omega)\subset L^2(\Omega)$ and spectral theory for self adjoint operators produce the desired result.

We point out that the eigenvalues of \eqref{eq:eigenvalue-original} can be expressed in terms of some explicit algebraic equations and the eigenfunctions admit an explicit representation in terms of the their respective eigenvalues.
For more details see the statement and the proof of \cite[Theorem 3.3]{Suspension} where the reader can realize that the eigenvalues can be classified into four different kinds.

\begin{figure}[th]
\begin{center}
 {\includegraphics[scale=0.20]{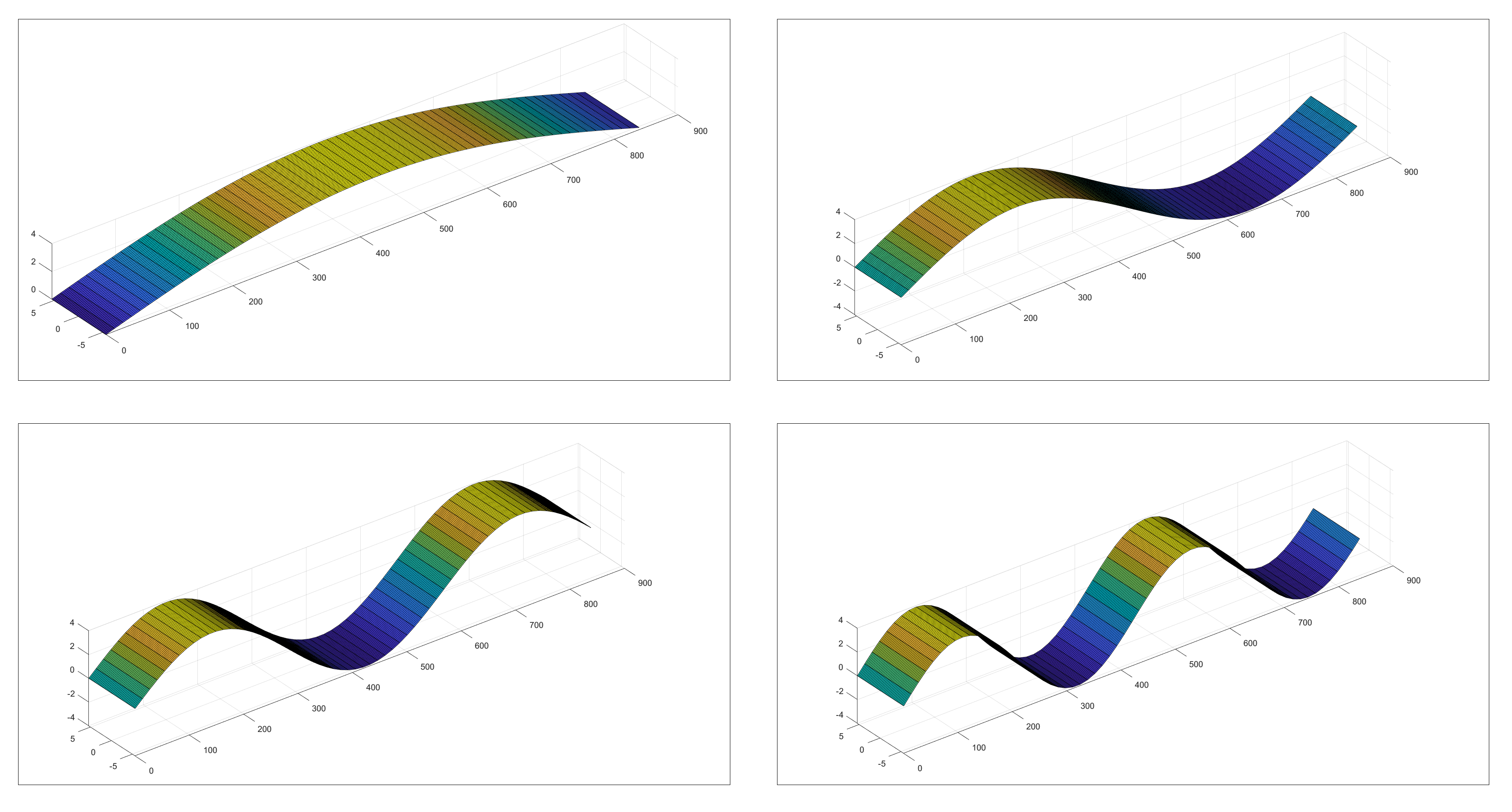}}
\caption{Eigenfunctions of $\lambda_1^{{\rm vert}}$, $\lambda_2^{{\rm vert}}$, $\lambda_3^{{\rm vert}}$,
$\lambda_4^{{\rm vert}}$.} \label{f:1}
\end{center}
\end{figure}

\begin{figure}[th]
\begin{center}
{\includegraphics[scale=0.20]{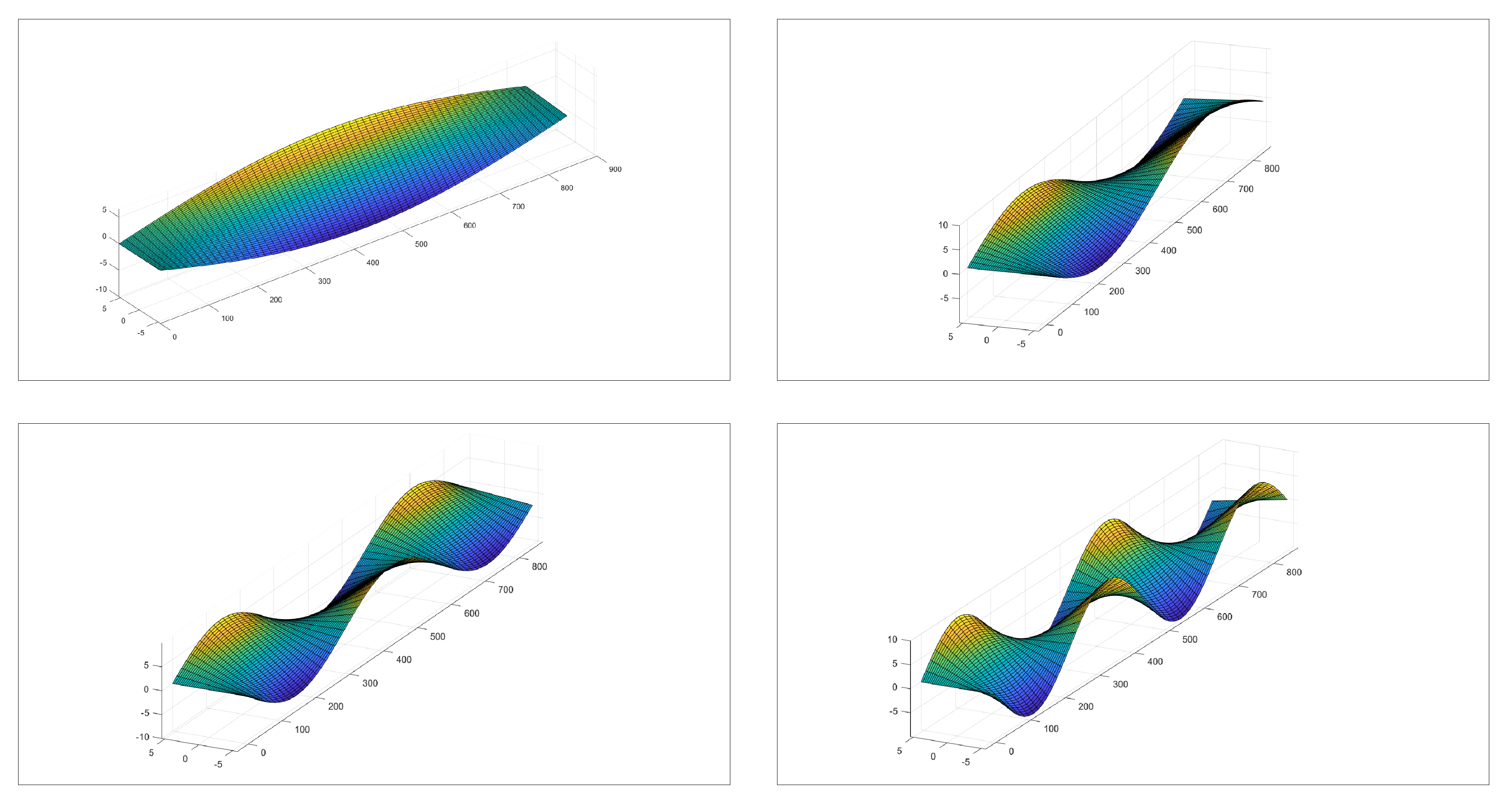}}
\caption{Eigenfunctions of $\lambda_1^{{\rm tors}}$, $\lambda_2^{{\rm tors}}$, $\lambda_3^{{\rm tors}}$,
$\lambda_4^{{\rm tors}}$.} \label{f:2}
\end{center}
\end{figure}

In \cite[Section 6]{Suspension}, we selected two families of eigenvalues denoted there by $\{\lambda_m^{{\rm vert}}\}_{m\ge 1}$ and $\{\lambda_m^{{\rm tors}}\}_{m\ge 1}$ respectively. Numerical evidence has shown that for the specific values assigned to the parameters of the plate, see \eqref{eq:elastic-par} below, the first eighteen eigenvalues all lie in one of these two families. We clarify that for any $m\ge 1$, the eigenfunctions corresponding to $\lambda_m^{{\rm vert}}$ are even with respect to the $y$ variable and the ones corresponding to $\lambda_m^{{\rm tors}}$ are odd with respect to the $y$ variable, thus giving sense to that notation (see Figure \ref{f:1} and Figure \ref{f:2}).

Let us consider the equation of motion for a free orthotropic plate:

\begin{equation} \label{eq:evolution}
  \frac{M}{2\ell} \, \frac{\partial^2 u}{\partial t^2}+\frac{d^3 \mathcal K}{12}
   \left(\Delta^2 u+\kappa \, \frac{\partial^4 u}{\partial x^4}\right)=0
\end{equation}
where $M$ is the mass linear density of the deck as explained in the introduction.

Given an eigenfunction $U_\lambda$ of \eqref{eq:eigenvalue-original} corresponding to some eigenvalue $\lambda$ we can construct a stationary wave solution in the form
\begin{equation} \label{eq:stat-wave}
  u_\lambda(x,y,t)=\sin(\omega_\lambda t) \, U_\lambda(x,y)
\end{equation}
where $\omega_\lambda$ represents an angular velocity. The frequency $\nu_\lambda$ is then obtained from $\omega_\lambda$ by dividing it by $2\pi$.

Inserting \eqref{eq:stat-wave} into \eqref{eq:evolution} and exploiting the fact that $U_\lambda$ is an eigenfunction with eigenvalue $\lambda$, we obtain
\begin{equation*}
   \left(-\frac{M}{2\ell} \, \omega_\lambda^2 +\lambda\right) \sin(\omega_\lambda t) U_\lambda(x,y)=0
\end{equation*}
and, in turn,
\begin{equation} \label{eq:nu-m}
  \nu_\lambda=\frac{\omega_\lambda}{2\pi}=\frac{1}{\pi} \sqrt{\frac{\ell \lambda}{2M}} \, .
\end{equation}

Looking at \eqref{eq:nu-m} and the definition of $\lambda_m^{{\rm vert}}$ and $\lambda_m^{{\rm tors}}$, it appears natural to define the following frequencies of vertical and torsional vibration respectively by
\begin{equation*}
  \nu_m^{{\rm vert}}=\frac{1}{\pi} \sqrt{\frac{\ell \lambda_m^{{\rm vert}}}{2M}} \, ,
  \qquad  \nu_m^{{\rm tors}}=\frac{1}{\pi} \sqrt{\frac{\ell \lambda_m^{{\rm tors}}}{2M}} \, .
\end{equation*}

Having in mind the structure of the Tacoma Narrows Bridge built and collapsed in 1940, taking inspiration from \cite{ammann, AG15, AG17}, in \cite{Suspension} we assigned to the parameters of the plate the following values
\begin{align}  \label{eq:elastic-par}
   & L=853.44 \ m \, , \qquad \ell=6 \ m \, , \qquad M=7198 \ kg/m \, , \qquad  \nu=0.2 \, , \\[7pt]
  \notag & E_1=2.1 \cdot 10^{11} \ Pa \, , \qquad E_2=1.687\cdot 10^9 \ Pa \, , \qquad  \mathcal R=2.109 \cdot 10^7 \ N\cdot m  \, , \qquad \kappa=123.48 \, ,
\end{align}
where $\mathcal R=\frac{d^3 \mathcal K}{12}$ denotes the rigidity of the plate.

Recalling equation \eqref{eq:familiar} and looking at \eqref{eq:elastic-par}, the reader can realize that the model of plate obtained above, exhibits a strongly anisotropic behavior since the value of the Young modulus $E_1$ is of two orders of magnitude larger than the Young modulus $E_2$. This suggests that the more popular isotropic plate model is not completely suitable to describe torsional oscillations of the deck of a bridge.

Assuming \eqref{eq:elastic-par} and exploiting the implicit representation of the eigenvalues provided in \cite{Suspension}, we numerically computed the values of the first ten frequencies of vertical oscillation and the first eight frequencies of torsional oscillation. Such values are collected in Table \ref{t:1}.

We observe that these frequencies are clearly relatively smaller than the ones expected for the oscillation of the deck of a suspension bridge. This fact is not strange at all since the oscillations that can be observed in a real suspension bridge are absolutely affected by the dynamics of cables and hangers.

This suggests that the present article and the related paper \cite{Suspension} have to be considered only as preliminary works
in the formulation of a complete model of suspension bridge in which the deck is described by an orthotropic plate.
A general idea in this direction was given in \cite{Suspension} where the following system was proposed:
\begin{equation} \label{eq:sistemone}
  \begin{cases}
     m \xi(x) \frac{\partial^2 p_1}{\partial t^2}-\frac{H_0}{(\xi(x))^2} \, \frac{\partial^2 p_1}{\partial x^2}
     =f_1\left(x,p_1,\frac{\partial p_1}{\partial x}\right)+F(u(\cdot,\ell)-p_1)  \\[7pt]
      m \xi(x) \frac{\partial^2 p_2}{\partial t^2}-\frac{H_0}{(\xi(x))^2} \, \frac{\partial^2 p_2}{\partial x^2}
     =f_2\left(x,p_2,\frac{\partial p_2}{\partial x}\right)+F(u(\cdot,-\ell)-p_2)  \\[7pt]
     \frac{M}{2\ell} \frac{\partial^2 u}{\partial t^2}+\frac{d^3 \mathcal K}{12}\left(\Delta^2 u+\kappa
     \frac{\partial^4 u}{\partial x^4}\right)=-F(u(\cdot,\ell)-p_1)-F(u(\cdot,-\ell)-p_2) \, .
  \end{cases}
\end{equation}
In \eqref{eq:sistemone}, $p_1=p_1(x,t)$ and $p_2(x,t)$ describes the displacements of the two cables, $m$ is the mass linear density of the cables, $s=s(x)$ is the configuration of the cables at rest, $\xi(x)=\sqrt{1+(s'(x))^2}$ is the local length of cables at rest, $H_0$ is the horizontal component of the tension of cables and $f_1, f_2, F$ are suitable nonlinearities to be determined in dependence of the accuracy one aims to achieve in the model.

The formulation of \eqref{eq:sistemone} was inspired by the models obtained in \cite{AG15,AG17}.

\begin{table}
  \begin{center}
    \begin{tabular}{|c|c|c|c|c|c|c|c|c|c|c|}
      \hline
         $m$                    & $1$       & $2$       & $3$      & $4$       & $5$       & $6$      & $7$ & $8$ & $9$ & $10$ \\
      \hline
         $\nu_m^{{\rm vert}}$   & $0.0045$  & $0.0180$  & $0.0406$ & $0.0722$  & $0.1128$ & $0.1624$ & $0.2211$ &
         $0.2887$ & $0.3654$    & $0.4512$ \\
      \hline
         $\nu_m^{{\rm tors}}$   & $0.0404$  & $0.0822$  & $0.1270$ & $0.1760$  & $0.2301$  & $0.2904$
                                & $0.3574$  & $0.4317$  & $-$ & $-$     \\
      \hline
    \end{tabular}
    \bigskip
  \caption{First ten frequencies corresponding to vertical oscillations measured in $Hz$}  \label{t:1}
 \end{center}
\end{table}


\bigskip

{\bf Acknowledgments} The author is member of the Gruppo Nazionale per l'Analisi Matematica, la Probabilit\`{a} e le loro Applicazioni (GNAMPA) of the Istituto Nazionale di Alta Matematica (INdAM). The author acknowledges partial financial support from the PRIN project 2017 ``Direct and inverse problems for partial differential equations: theoretical aspects and
applications'' and from the INDAM - GNAMPA project 2019 ``Analisi spettrale per operatori ellittici del secondo e quarto ordine con condizioni al contorno di tipo Steklov o di tipo parzialmente incernierato''.

This research was partially supported by the research project ``Metodi e modelli per la matematica e le sue
applicazioni alle scienze, alla tecnologia e alla formazione'' Progetto di Ateneo 2019 of the University of Piemonte Orientale ``Amedeo Avogadro''.

The author is grateful to Elvise Berchio, Alessio Falocchi and Pier Domenico Lamberti for the useful discussions and suggestions that supported the beginning of this work.

\end{document}